\documentclass[runningheads,a4paper]{llncs}
\usepackage{amssymb}
\usepackage{graphicx}
\usepackage{url}
\newcommand{\keywords}[1]{\par\addvspace\baselineskip
\noindent\keywordname\enspace\ignorespaces#1}

\usepackage{tikz}
\usepackage{pgfplots}
\usepackage{amsmath}
\definecolor{veryverylight}{gray}{0.85}
\definecolor{verylight}{gray}{0.75}
\definecolor{light}{gray}{0.65}
\definecolor{medium}{gray}{0.55}

\usepackage{paralist}

\usepackage{algorithm,caption} 
\usepackage[noend]{algpseudocode}

\algrenewcommand\algorithmicrequire{\textbf{Input:}}
\algrenewcommand\algorithmicensure{\textbf{Return:}}

\def\ZZ{{\mathbb Z}}

\def\RR{{\mathbb R}}

\def\cone{\operatorname{cone}}

\def\Hilb{\operatorname{Hilb}}

\def\gp{\operatorname{gp}}

\def\aff{\operatorname{aff}}

\begin{document}

\mainmatter

\title{The subdivision of large simplicial cones in Normaliz}  
\titlerunning{Recent developments in Normaliz} 
\author{Winfried Bruns\inst{1} \and  Richard Sieg\inst{2} \and Christof S\"oger\inst{3}}
\authorrunning{Bruns-Sieg-S\"oger}
\institute{
University of Osnabr\"uck, Germany\\
\email{wbruns@uos.de},\\
\texttt{http://www.home.uni-osnabrueck.de/wbruns/}
\and
University of Osnabr\"uck, Germany\\
\email{risieg@uos.de},\\
\texttt{http://www.math.uni-osnabrueck.de/normaliz/}
\and
University of Osnabr\"uck, Germany\\
\email{csoeger@uos.de},\\
\texttt{http://www.math.uni-osnabrueck.de/normaliz/}
}
\maketitle

\begin{abstract}
Normaliz is an open-source software for the computation of lattice points in rational polyhedra, or, in a different language, the solutions of linear diophantine systems. The two main computational goals are (i) finding a system of generators of the set of lattice points and (ii) counting elements degree-wise in a generating function, the Hilbert Series. In the homogeneous case, in which the polyhedron is a cone, the set of generators is the Hilbert basis of the intersection of the cone and the lattice, an affine monoid.

We will present some improvements to the Normaliz algorithm by subdividing simplicial cones with huge volumes. In the first approach the subdivision points are found by integer programming techniques. For this purpose we interface to the integer programming solver SCIP to our software. In the second approach we try to find good subdivision  points in an approximating overcone that is faster to compute. 

\end{abstract}


\section{Introduction}

Normaliz \cite{Normaliz} is a software for the computation of lattice points in rational polyhedra. These are exactly the solutions of linear diophantine systems of inequalities, equations and congruences. It pursues two main computational goals:
(i) finding a minimal generating system of the set of lattice points in a polyhedron; (ii) counting elements degree-wise in a generating function, the Hilbert series. In the homogeneous case, in which the polyhedron is a cone, the set of generators is the Hilbert basis of the intersection of the cone and the lattice, which is an affine monoid by Gordan's lemma. For the mathematical background we refer the reader to \cite{BGu}. The Normaliz algorithms are described in \cite{BI} and \cite{BIS}. The second paper contains extensive performance data.

Normaliz (present public version 3.1.1) is written in C++ (using Boost and GMP/MPIR), parallelized with OpenMP, and runs under Linux, MacOs and MS Windows. It is based on its C++ library libnormaliz which offers the full functionality of Normaliz.
There are file based interfaces for Singular, Macaulay 2 and Sage, and C++ level interfaces  for CoCoA, polymake, Regina and GAP. A C++ level interface to Sage should be available in the near future. There is also the GUI interface jNormaliz.

Normaliz has found applications in commutative algebra, toric geometry,
combinatorics, integer programming, invariant theory, elimination
theory, group theory, mathematical logic, algebraic topology and even theoretical physics.

\section{Hilbert basis and Hilbert series}

We will first describe the main functionality of Normaliz. For simplicity we restrict ourselves to homogeneous linear systems in the following, or, geometrically speaking, to the intersections of lattices $L\subset \ZZ^d$ and rational cones $C\subset \RR^d$. 

\begin{definition}
A (rational) \emph{polyhedron} $P$ is the intersection of finitely many (rational) halfspaces. If it is bounded, then it is called a \emph{polytope}. If all the halfspaces are linear, then $P$ is a cone.

The \emph{dimension} of $P$ is the dimension of the smallest affine subspace $\aff(P)$ containing $P$.

An \emph{affine monoid} is a finitely generated submonoid of $\ZZ^d$ for some $d$.
\end{definition}
By the theorem of Minkowski-Weyl, $C\subset\RR^d$ is a (rational) cone if and only if there exist finitely many (rational) vectors $x_1,\dots,x_n$ such that
$$
C=\operatorname{cone}(x_1,\dots,x_n)=\{a_1x_1+\dots+a_nx_n:a_1,\dots,a_n\in\RR_+\}.
$$
If $x_1,\ldots,x_n$ are linearly independent, we call $C$ \emph{simplicial}. For Normaliz, cones $C$ and lattices $L$ can either be specified by generators $x_1, \dots, x_n \in \ZZ^d$ or by constraints, i.e., homogeneous systems of diophantine linear inequalities, equations and congruences.
Normaliz also offers to define an affine monoid as the quotient  of $\ZZ_+^n$ modulo the  intersection with a sublattice of $\ZZ^n$. 

Normaliz puts no restriction on the rational cone $C$. In the following we will however assume that $C$ is pointed, i.e. $x,-x\in C \Rightarrow x=0$. This is justified since computations in non-pointed cones are done via the projection to the quotient modulo the maximal linear subspace, which is pointed.

\begin{figure}[hbt]
\begin{center}

\includegraphics{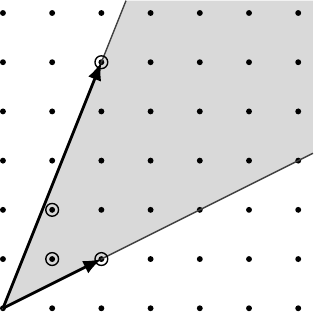}$\quad$\includegraphics{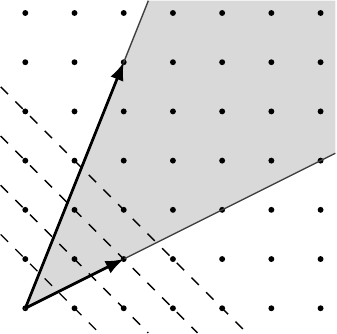}

\end{center}
\caption{A cone with the Hilbert basis (circled points) and grading.}\label{figCone}
\end{figure}

By Gordan's lemma the monoid $M = C \cap L$ is finitely generated. This affine monoid has a (unique) minimal generating system called the \emph{Hilbert basis} $\Hilb(M)$, see Figure \ref{figCone} for an example. The computation of the Hilbert basis is the first main task of Normaliz.

One application is the computation of the \emph{normalization} of an affine monoid $M$; this explains the name Normaliz. The normalization is the intersection of the cone generated by $M$ with the sublattice $\gp(M)$ generated by $M$. One calls $M$ \emph{normal}, if it coincides with its normalization.

The second main task is to compute the Hilbert (or Ehrhart) series of a graded monoid. A \emph{grading} of a
monoid $M$ is simply a homomorphism $\deg:M\to\ZZ^g$ where
$\ZZ^g$ contains the degrees. The \emph{Hilbert series} of $M$
with respect to the grading is the formal Laurent series
$$
H(t)=\sum_{u\in \ZZ^g} \#\{x\in M: \deg x=u\}t_1^{u_1}\cdots t_g^{u_g}=\sum_{x\in M}t^{\deg x},
$$
provided all sets $\{x\in M: \deg x=u\}$ are finite. At the moment, Normaliz can only handle the case $g=1$, and therefore we restrict ourselves to this case.
We assume in the following that $\deg x >0$ for all nonzero $x\in M$ and that there exists an $x\in\gp(M)$ such that $\deg x=1$. (Normaliz always rescales the grading accordingly.)

Assume that $M$ is a normal and affine monoid. By a theorem of Hilbert and Serre (\cite[Theorem 6.37]{BGu}), $H(t)$ in the $\ZZ$-graded case 
is the Laurent expansion of a rational function at the origin:
$$
H(t)=\frac{R(t)}{(1-t^e)^r},\qquad R(t)\in\ZZ[t], 
$$
where $r$ is the rank of $M$ and $e$ is the least common multiple
of the degrees of the extreme integral generators of $\cone(M)$. As a rational function, $H(t)$ has negative degree.


A rational cone $C$ and a grading together define the rational
polytope $Q=C\cap A_1$ where $A_1=\{x:\deg x=1\}$. In this
sense the Hilbert series is nothing but the Ehrhart series of
$Q$.



\section{The primal algorithm}\label{sec:algo}

The \emph{primal} Normaliz algorithm is triangulation based. Normaliz contains a second, \emph{dual} algorithm for the computation of Hilbert bases that implements ideas of Pottier \cite{Pot}. The dual algorithm is treated in \cite{BI}, and has not changed much in the last years and we do not discuss it in this article.

The primal algorithm starts from a pointed rational cone $C\subset\RR^d$ given by a system of generators $x_1,\dots,x_n$ and a sublattice $L\subset\ZZ^d$ that contains $x_1,\dots,x_n$. Other types of input data are first transformed into this format. The algorithm is composed as follows:
\begin{enumerate}
\item Initial coordinate transformation to $E=L\cap (\RR x_1+\dots+\RR x_n)$;
\item Fourier-Motzkin elimination computing the support hyperplanes of $C$;
\item computation of a triangulation, i.e. a face-to-face decomposition into simplicial cones;
\item evaluation of the simplicial cones in the triangulation;
\item collection of the local data;
\item reverse coordinate transformation to $\ZZ^d$.
\end{enumerate}

The algorithm does not strictly follow this chronological order, but interleaves steps 2--5 in an intricate way to ensure low memory usage and efficient parallelization. 

%


\subsection{Simplicial cones}

We will now focus on step 4 of the primal algorithm, the evaluation of simplicial cones. Let $x_1,\dots,x_d\in\ZZ^d$ be linearly independent and $S=\operatorname{cone}(x_1,\dots,x_d)$.
Then the integer points in the \emph{fundamental domain} of $S$ 
$$
E=\{q_1x_1+\dots+q_dx_d: 0\le q_i<1\}\cap\ZZ^d
$$
together with $x_1,\dots,x_d$ generate the monoid $S\cap\ZZ^d$. 
\begin{figure}[h]
\begin{center}
\includegraphics{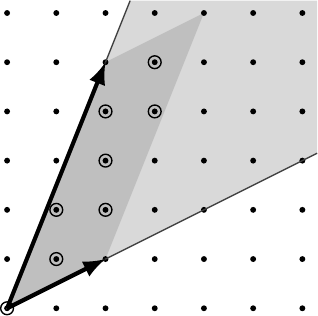}
\end{center}
\caption{A cone with a fundamental domain}
\end{figure}

Every residue class in the quotient $\ZZ^d/U$, where $U=\ZZ x_1+\dots+\ZZ x_d$, has exactly one representative in $E$.
Representatives of residue classes can be quickly computed via the elementary divisor algorithm and
from an arbitrary representative we obtain the one in $E$ by division with remainder. The integer points of the fundamental domain are candidates for the Hilbert basis of the cone. After their computation they are shrunk to the Hilbert basis by successively discarding elements $x$ which are \emph{reducible}, i.e. there exists an $y\in E, y\neq x$ such that $x-y\in C$. Also the computation of the Hilbert series uses the set $E$ and a \emph{Stanley decomposition} based on it; see \cite{BIS}.

The number of elements in $E$ is given by the (lattice normalized) volume of the simplex:
$$ 
|E| = \operatorname{vol}(S) = \det(x_1,\dots,x_d).
$$
Therefore the determinant of the generators of the simplicial cone has an enormous impact on the runtime of the Normaliz algorithm. The algorithms presented in this paper try to decompose a simplex with big volume into simplices such that the sum of their volumes is considerably smaller. For this purpose we compute integer points from the cone and use them for a new triangulation. 

Theoretically the best choice for these points are the vertices of the \emph{bottom} $B(S)$ of the simplex which is defined as the union of the bounded faces of the polyhedron  $\operatorname{conv}((S\cap\ZZ^d)\setminus\{0\})$.
In practice, the computation of the whole bottom would equalize the benefit from the small volume or even make it worse. 

Therefore, we determine only some points from the bottom. Normaliz employs two methods for this purpose:
\begin{compactenum}[(1)]
\item computation of subdivision points by integer programming methods,
\item computation of candidate subdivision points by approximation of the given simplicial cone by an overcone that is generated by vectors  of ``low denominator''.
\end{compactenum}

\section{Methods from integer programming}\label{sec:ip}

\begin{figure}[hbt]
\includegraphics{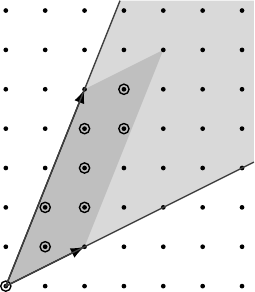}
     \hspace{0.3cm}
     \includegraphics{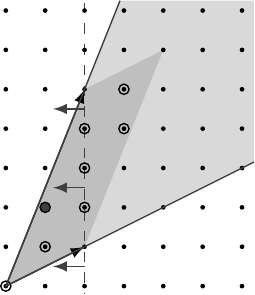}
     \hspace{0.3cm}
     \includegraphics{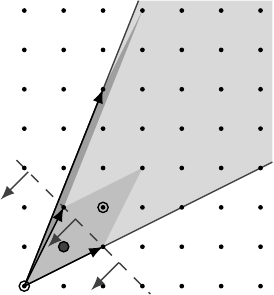}
     \hspace{0.3cm}
	\includegraphics{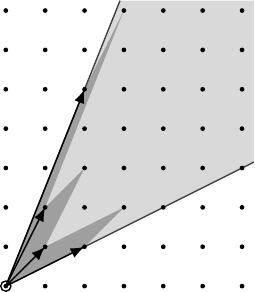}
     \caption{The integer programming algorithm for a cone}\label{fig:ipalgo}
\end{figure}

For each simplex $S=\operatorname{cone}(x_1,\ldots,x_d)$ in the triangulation with large enough volume we try to compute a point $x$ that minimizes the \emph{sum of determinants}: 
$$
\sum_{i=1}^{d}\det(x_1,\ldots,x_{i-1},x,x_{i+1},\ldots,x_d),
$$
which can also be expressed as $N^T x$, where $N$ is a normal vector on the affine hyperplane spanned by $x_1,\ldots,x_d$. Such a  point can be found by solving the following integer program:
$$
\qquad\min\{N^T x : x\in S\cap\mathbb{Z}^d,x\neq 0,N^T x< N^T x_1\}.\qquad(\star)
$$     
If the problem has a solution $\hat{x}$, we form a  \emph{stellar subdivision} of the simplex with respect to $\hat{x}$: For every support hyperplane $H_i$ (not containing $x_i$) which does not contain $\hat{x}$ we form the simplex
$$
T_i=\operatorname{cone}(x_1,\ldots,x_{i-1},\hat{x},x_{i+1},\ldots,x_d).
$$
If the volume of $T_i$ is larger than a particular bound, we repeat this process and continue until all simplices have a smaller volume than this bound or the corresponding integer problems have no solutions. Figure~\ref{fig:ipalgo} illustrates the algorithm.

After computing a set of integer points $\mathcal{B}$, we triangulate the bottom of $\operatorname{conv}(\mathcal{B}\cup\{x_1,\ldots,x_d\})$ and continue by evaluating this triangulation with the usual Normaliz algorithm.

\subsection{Implementation and results}
We use the mixed integer programming solver SCIP \cite{Scip} via its C++ interface. The algorithm runs in parallel with one SCIP environment for every thread using OpenMP. Moreover each SCIP instance has its own time limit ($\log(\operatorname{vol}(S))^2$ sec) and feasibility bounds. 

The condition that $x\neq 0$ could be implemented by the inequality $N^T x\geq1$. However this approach is prone to large numbers in $N$. Therefore we first check, whether all generators are positive in one entry $i$ and thus require $x_i\geq1$. If this is not the case we make a bound disjunction of the form $(x_i\leq -1 \vee x_i\geq 1)$.

Table~\ref{tab:ip} presents example data computed on a SUN xFire 4450 with four Intel Xeon X7460 processors, using 20 threads and solving integer programs only for simplices with a volume larger than $10^6$.
\captionsetup[table]{skip=8pt}
\begin{table}[h]

\centering
\setlength{\tabcolsep}{3.2pt}
\renewcommand{\arraystretch}{1.2}
\begin{tabular}{|c|c|c|c|}
\hline
 & hickerson-16 & hickerson-18 & knapsack\_11\_60  \\ \hline
simplex volume & $9.83\times 10^7$ & $4.17\times 10^{14}$ & $2.8\times 10^{14}$ \\ \hline
volume under bottom & {$8.10\times 10^5$} & {$3.86\times 10^7$} & {$2.02\times 10^7$} \\ \hline 
volume used  & {$3.93\times 10^6$} & {$5.47\times 10^7$} &  {$2.39\times 10^7$} \\ \hline
integer programs solved & $4$ & $582016$ & $11621$ \\ \hline
improvement factor & {25} &  {$7.62\times10^6$} & {$1.17\times 10^7$}\\ \hline
runtime without subdivision &  {2s} & {$>12$d} &  {$>8$d} \\ \hline
runtime with subdivision  &  {0.5s} & {46s} & {5.1s} \\ \hline
\end{tabular}
\caption{Runtime improvements using integer programming methods}\label{tab:ip}
\end{table}

\begin{figure}[h]
\centering
\includegraphics{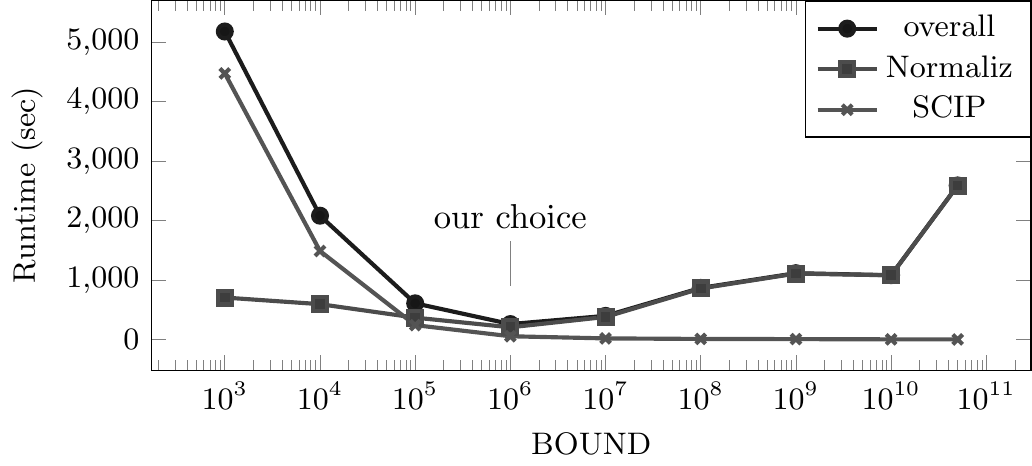}
\caption{Runtime graph showing different choices for the bound}\label{fig:bound}
\end{figure}

The bound on the volume to stop the calculation of a single simplex has a significant effect on the runtime of the algorithm. A smaller bound means that more integer programs have to be solved by SCIP, whereas a large bound prevents a major improvement of the respective volume. Running several experiments, it turns out that $10^6$ is a good value in between these two extreme cases. Figure~\ref{fig:bound} shows a runtime graph illustrating the effect of different choices for this bound. The measured time is a single thread computation of hickerson-18.

\section{Approximation}
SCIP cannot be employed in all environments. Especially if Normaliz is bundled with another software package it may be undesirable or even impossible to force the link to SCIP.  

Our second approach is completely implemented within Normaliz. It first approximates the simplicial cone $S$ by a (not necessarily simplicial) overcone $C$ for which the sets $E$ in a triangulation of $C$ are significantly faster to compute. Then these points are used to decompose the original simplex as before. It is clear that the efficiency depends crucially on the intersection of the sets $E$ with $S$. 

\begin{figure}[h]
\begin{center}

\includegraphics{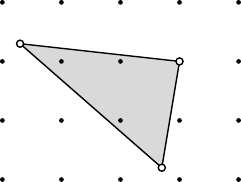}
\hspace{0.7cm}
\includegraphics{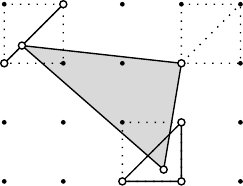}
\hspace{0.7cm}
\includegraphics{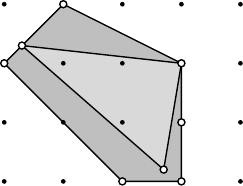}
   \end{center}
\caption{Approximating cone}\label{fig:approx}
\end{figure}

For this purpose we look at the polytope given by the cross section of the simplex at height one, where the height function comes from the normal vector $N$ on the affine hyperplane spanned by the generators. For every vertex of this polytope we triangulate the lattice cube around it using the braid hyperplane arrangement $\{x_i=x_j\}$. We continue by detecting the minimal face  containing the vertex and collect its vertices, which are at most $d$. The approximating cone $C$ is then generated by all vertices found in that way. Figure~\ref{fig:approx} illustrates the choice of the approximation for a 3-dimensional cone (with a 2-dimensional cross section).

\begin{figure}[h]
\begin{center}
\includegraphics{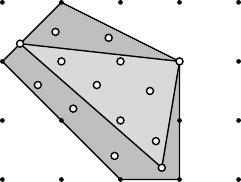}
\hspace{0.7cm}
\includegraphics{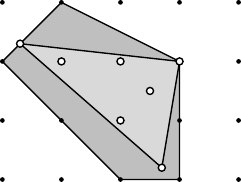}
\hspace{0.7cm}
\includegraphics{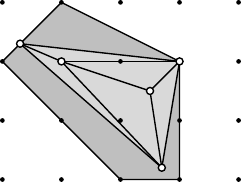}
   \end{center}
\caption{Decomposition of a simplex after approximation}\label{fig:approx_cand}
\end{figure}

As in the usual Normaliz algorithm we create a candidate list for the exterior cone, but keep only those points which lie inside the original simplex $S$. The remaining candidates are then reduced as before, which results in a list $\mathcal{B}$ which is used for a recursive decomposition of the simplex as in Section~\ref{sec:ip}. Figure~\ref{fig:approx_cand} illustrates this process for the previous example.

It might happen for both algorithms that no decomposition point can be found, although the volume of the simplex is still quite large ($>10^9$) and subdivision points exist. In this case, the approximation method is applied again with a higher level of approximation. 


Table~\ref{tab:approx} contains performance data for the examples in Section~\ref{sec:ip}.
\begin{table}[h]
\centering
\setlength{\tabcolsep}{3.2pt}
\renewcommand{\arraystretch}{1.2}
\begin{tabular}{|c|c|c|c|}
\hline
 & hickerson-16 & hickerson-18 & knapsack\_11\_60  \\ \hline
 volume used  & {$3.93\times 10^6$} & {$8.42\times 10^7$} &  {$9.36\times 10^7$} \\ \hline
 improvement factor & {25} &  {$4.95\times10^6$} & {$2.99\times 10^4$}\\ \hline
runtime with subdivision  &  {0.4s} & {50s} & {2m30s} \\ \hline
\end{tabular}
\caption{Runtime improvements using the approximation method}\label{tab:approx}
\end{table}
\vspace{-1cm}

At present we are working on improvements of the approximation method.


\end{document}